\documentclass{amsart}
\usepackage{amstext}
\usepackage{amsthm}
\usepackage{amssymb}
\usepackage{esint}

\makeatletter

\numberwithin{equation}{section}
\numberwithin{figure}{section}
\theoremstyle{plain}
\newtheorem{thm}{Theorem}

\allowdisplaybreaks
\usepackage[sort&compress,square,numbers]{natbib}
\usepackage{tikz}
\usetikzlibrary{matrix,positioning}

\makeatother

\begin{document}

\title[Identities involving Bessel polynomials]{Identities involving Bessel polynomials arising from linear differential
equations}

\author{Taekyun Kim}
\address{Department of Mathematics, Kwangwoon University, Seoul 139-701, Republic
of Korea}
\email{tkkim@kw.ac.kr}

\author{Dae San Kim}
\address{Department of Mathematics, Sogang University, Seoul 121-742, Republic
of Korea}
\email{dskim@sogang.ac.kr}

\begin{abstract}
In this paper, we study linear differential equations arising from Bessel
polynomials and their applications. From these linear differential
equations, we give some new and explicit identities for Bessel polynomials.
\end{abstract}


\maketitle
\global\long\def\relphantom#1{\mathrel{\phantom{{#1}}}}

\section{Introduction}

As is well known, the Bessel differential equation is given by
\begin{equation}
x^{2}\frac{d^{2}y}{dx^{2}}+x\frac{dy}{dx}+\left(x^{2}-\alpha^{2}\right)y=0,\quad\left(\text{see \cite{key-17}}\right).\label{eq:1-1}
\end{equation}
for an arbitrary complex number $\alpha$.

The Bessel functions of the first kind $J_{\alpha}\left(x\right)$
are defined by the solution of (\ref{eq:1-1}).

For $n\in\mathbb{Z}$, $J_{n}\left(x\right)$ are sometimes also called
cylinder function or cylindrical harmonics.

It is known that
\begin{equation}
J_{n}\left(x\right)=\sum_{l=0}^{\infty}\frac{\left(-1\right)^{l}}{l!\left(n+l\right)!}\left(\frac{x}{2}\right)^{2l+n},\quad\left(\text{see \cite{key-15,key-16,key-17}}\right).\label{eq:2-1}
\end{equation}

The generating function of Bessel functions is given by
\begin{equation}
e^{\frac{x}{2}\left(t-\frac{1}{t}\right)}=\sum_{n=-\infty}^{\infty}J_{n}\left(x\right)t^{n},\label{eq:3-1}
\end{equation}
and $J_{n}\left(x\right)$ can be also represented by the contour
integral as
\begin{equation}
J_{n}\left(x\right)=\frac{1}{2\pi i}\oint e^{\frac{x}{2}\left(t-\frac{1}{t}\right)}t^{-n-1}dt,\quad\left(\text{see \cite{key-17}}\right),\label{eq:4-1}
\end{equation}
where the contour encloses the origin and is traversed in a counterclockwise
direction.

The Bessel polynomials are defined by the solution of the differential
equation
\begin{equation}
x^{2}\frac{d^{2}y}{dx^{2}}+2\left(x+1\right)\frac{dy}{dx}-n\left(n+1\right)y=0,\quad\left(\text{see \cite{key-1,key-2,key-3,key-4,key-5,key-14,key-15,key-16}}\right).\label{eq:5-1}
\end{equation}

Indeed, the solutions of (\ref{eq:5-1}) are given by
\begin{align}
y_{n}\left(x\right) & =\sum_{k=0}^{n}\frac{\left(n+k\right)!}{\left(n-k\right)!k!}\left(\frac{x}{2}\right)^{k}\label{eq:6-1}\\
 & =\sqrt{\frac{2}{\pi x}}e^{\frac{1}{x}}K_{-n-\frac{1}{2}}\left(\frac{1}{x}\right),\quad\left(\text{see \cite{key-14,key-15,key-16,key-17}}\right),\nonumber
\end{align}
where
\[
K_{\nu}\left(z\right)=\frac{\Gamma\left(\nu+\frac{1}{2}\right)\left(2z\right)^{\nu}}{\sqrt{\pi}}\int_{0}^{\infty}\frac{\cos t}{\left(t^{2}+z^{2}\right)^{\nu+\frac{1}{2}}}dt.
\]

We note that $y_{n}\left(x\right)$ are very similar to the modified
spherical Bessel function of the second kind.

The first few are given as
\begin{align*}
y_{0}\left(x\right) & =1,\quad y_{1}\left(x\right)=x+1,\quad y_{2}\left(x\right)=3x^{2}+3x+1,\\
y_{3}\left(x\right) & =15x^{3}+15x^{2}+6x+1,\\
y_{4}\left(x\right) & =105x^{4}+105x^{3}+45x^{2}+10x+1,\quad\dots.
\end{align*}

Carlitz reverse Bessel polynomials are defined by
\begin{equation}
p_{n}\left(x\right)=x^{n}y_{n-1}\left(\frac{1}{x}\right),\quad\left(n\in\mathbb{N}\cup\left\{ 0\right\} \right),\quad\left(\text{see \cite{key-3,key-14}}\right).\label{eq:7-1}
\end{equation}

These polynomials are also given by the generating function as
\begin{equation}
e^{x\left(1-\sqrt{1-2t}\right)}=\sum_{n=0}^{\infty}p_{n}\left(x\right)\frac{t^{n}}{n!}.\label{eq:8-1}
\end{equation}

The explicit formulas for them are
\begin{align}
p_{n}\left(x\right) & =\sum_{k=1}^{n}\frac{\left(2n-k-1\right)!}{2^{n-k}\left(k-1\right)!\left(n-k\right)!}x^{k}\label{eq:9-1}\\
 & =\left(2n-3\right)!!x\,_{1}F_{1}\left(1-n;2-2n;2x\right),\quad\left(\text{see \cite{key-14,key-15,key-16}}\right),\nonumber
\end{align}
where
\[
n!!=\begin{cases}
n\left(n-2\right)\cdots5\cdot3\cdot1 & \text{if }n>0\text{ odd},\\
n\left(n-2\right)\cdots6\cdot4\cdot2 & \text{if }n>0\text{ even,}\\
1 & \text{if }n=-1,0,
\end{cases}
\]
and
\begin{align*}
\,_{1}F_{1}\left(a;b;z\right) & =1+\frac{a}{b}z+\frac{a\left(a+1\right)}{b\left(b+1\right)}\frac{z^{2}}{2!}+\cdots\\
 & =\sum_{k=0}^{\infty}\frac{a\left(a+1\right)\cdots\left(a+k-1\right)}{b\left(b+1\right)\cdots\left(b+k-1\right)}\frac{z^{k}}{k!}\\
 & =\frac{\Gamma\left(b\right)}{\Gamma\left(b-a\right)\Gamma\left(a\right)}\int_{0}^{1}e^{zt}t^{a-1}\left(1-t\right)^{b-a-1}dt.
\end{align*}

The first few polynomials are
\begin{align*}
p_{1}\left(x\right) & =x,\\
p_{2}\left(x\right) & =x^{2}+x,\\
p_{3}\left(x\right) & =x^{3}+3x^{2}+3x,\\
p_{4}\left(x\right) & =x^{4}+6x^{3}+15x^{2}+15x,\cdots.
\end{align*}

Recently, several authors have studied non-linear differential equations
related to special polynomials (see \cite{key-6,key-7,key-8,key-9,key-10,key-11,key-12,key-13}).

The reverse Bessel polynomials are used in the design of Bessel electronic
filters.

In this paper, we consider linear differential equations arising from
Carlitz reverse Bessel polynomials and give some new and explicit
identities for Bessel polynomials.

\section{Identities involving Bessel polynomials arising from linear differential
equations}

Let us put
\begin{equation}
F=F\left(t,x\right)=e^{x\left(1-\sqrt{1-2t}\right)}.\label{eq:1}
\end{equation}

Thus, by (\ref{eq:1}), we get
\begin{equation}
F^{\left(1\right)}=\frac{d}{dt}F\left(t,x\right)=x\left(1-2t\right)^{-\frac{1}{2}}F,\label{eq:2}
\end{equation}
\begin{align}
F^{\left(2\right)} & =\frac{dF^{\left(1\right)}}{dt}\label{eq:3}\\
 & =\left(x\left(1-2t\right)^{-\frac{3}{2}}+x^{2}\left(1-2t\right)^{-1}\right)F,\nonumber
\end{align}
\begin{align}
F^{\left(3\right)} & =\frac{d}{dt}F^{\left(2\right)}\label{eq:4}\\
 & =\left(3x\left(1-2t\right)^{-\frac{5}{2}}+3x^{2}\left(1-2t\right)^{-2}+x^{3}\left(1-2t\right)^{-\frac{3}{2}}\right)F,\nonumber
\end{align}
and
\begin{align}
F^{\left(4\right)} & =\frac{dF^{\left(3\right)}}{dt}\label{eq:5}\\
 & =\left(15x\left(1-2t\right)^{-\frac{7}{2}}+15x^{2}\left(1-2t\right)^{-3}+6x^{3}\left(1-2t\right)^{-\frac{5}{2}}+x^{4}\left(1-2t\right)^{-2}\right)F.\nonumber
\end{align}

Continuing this process, we set
\begin{align}
F^{\left(N\right)} & =\left(\frac{d}{dt}\right)^{N}F\left(t,x\right)\label{eq:6}\\
 & =\left(\sum_{i=N}^{2N-1}a_{i-N}\left(N,x\right)\left(1-2t\right)^{-\frac{i}{2}}\right)F,\nonumber
\end{align}
where $N=1,2,3,\dots$.

From (\ref{eq:6}), we note that
\begin{align}
 & F^{\left(N+1\right)}\label{eq:7}\\
 & =\frac{d}{dt}F^{\left(N\right)}\nonumber \\
 & =\left(\sum_{i=N}^{2N-1}a_{i-N}\left(N,x\right)\left(-\frac{i}{2}\right)\left(1-2t\right)^{-\frac{i}{2}-1}\left(-2\right)\right)F\nonumber \\
 & \relphantom =+\sum_{i=N}^{2N-1}a_{i-N}\left(N,x\right)\left(1-2t\right)^{-\frac{i}{2}}F^{\left(1\right)}\nonumber \\
 & =\left(\sum_{i=N}^{2N-1}ia_{i-N}\left(N,x\right)\left(1-2t\right)^{-\frac{i+2}{2}}\right)F\nonumber \\
 & \relphantom =+\left(\sum_{i=N}^{2N-1}a_{i-N}\left(N,x\right)\left(1-2t\right)^{-\frac{i}{2}}\right)x\left(1-2t\right)^{-\frac{1}{2}}F\nonumber \\
 & =\left(\sum_{i=N}^{2N-1}ia_{i-N}\left(N,x\right)\left(1-2t\right)^{-\frac{i+2}{2}}\right)F+\left(\sum_{i=N}^{2N-1}xa_{i-N}\left(N,x\right)\left(1-2t\right)^{-\frac{i+1}{2}}\right)F\nonumber \\
 & =\left\{ xa_{0}\left(N,x\right)\left(1-2t\right)^{-\frac{N+1}{2}}+\left(2N-1\right)a_{N-1}\left(N,x\right)\left(1-2t\right)^{-\frac{2N+1}{2}}\right.\nonumber \\
 & \relphantom =\left.+\sum_{i=N+1}^{2N-1}\left(\left(i-1\right)a_{i-N-1}\left(N,x\right)+xa_{i-N}\left(N,x\right)\right)\left(1-2t\right)^{-\frac{i+1}{2}}\right\} F.\nonumber
\end{align}

By replacing $N$ by $N+1$ in (\ref{eq:6}), we get
\begin{align}
F^{\left(N+1\right)} & =\left(\sum_{i=N+1}^{2N+1}a_{i-N-1}\left(N+1,x\right)\left(1-2t\right)^{-\frac{i}{2}}\right)F\label{eq:8}\\
 & =\left(\sum_{i=N}^{2N}a_{i-N}\left(N+1,x\right)\left(1-2t\right)^{-\frac{i+1}{2}}\right)F.\nonumber
\end{align}

By comparing the coefficients on both sides (\ref{eq:7}) and
(\ref{eq:8}), we have
\begin{align}
a_{0}\left(N+1,x\right) & =xa_{0}\left(N,x\right),\label{eq:9}\\
a_{N}\left(N+1,x\right) & =\left(2N-1\right)a_{N-1}\left(N,x\right),\label{eq:10}
\end{align}
and
\begin{equation}
a_{i-N}\left(N+1,x\right)=\left(i-1\right)a_{i-N-1}\left(N,x\right)+xa_{i-N}\left(N,x\right),\label{eq:11}
\end{equation}
where $N+1\le i\le2N-1$.

From (\ref{eq:2}) and (\ref{eq:6}), we can derive the following
equation (\ref{eq:11}):
\begin{equation}
x\left(1-2t\right)^{-\frac{1}{2}}F=F^{\left(1\right)}=a_{0}\left(1,x\right)\left(1-2t\right)^{-\frac{1}{2}}F.\label{eq:12}
\end{equation}

Thus, by (\ref{eq:12}), we have
\begin{equation}
a_{0}\left(1,x\right)=x.\label{eq:13}
\end{equation}

From (\ref{eq:9}), we note that
\begin{equation}
a_{0}\left(N+1,x\right)=xa_{0}\left(N,x\right)=x^{2}a_{0}\left(N-1,x\right)=\cdots=x^{N}a_{0}\left(1,x\right)=x^{N+1},\label{eq:14}
\end{equation}
and, by (\ref{eq:10}), we see
\begin{align}
a_{N}\left(N+1,x\right) & =\left(2N-1\right)a_{N-1}\left(N,x\right)\label{eq:15}\\
 & =\left(2N-1\right)\left(2N-3\right)a_{N-2}\left(N-1,x\right)\nonumber \\
 & \vdots\nonumber \\
 & =\left(2N-1\right)\left(2N-3\right)\cdots3\cdot1a_{0}\left(1,x\right)\nonumber \\
 & =\left(2N-1\right)!!x.\nonumber
\end{align}

The matrix $\left(a_{i}\left(j,x\right)\right)_{0\le i\le N-1,1\le j\le N}$
is given by

\begin{equation*}
 \begin{tikzpicture}[baseline=(current  bounding  box.west)]
  \matrix (mymatrix) [matrix of math nodes,left delimiter={[},right
delimiter={]}]
  {
  x & x^{2} & x^{3} & x^{4} & \cdots & x^{N}\\
 & 1!!x\\
 &  & 3!!x\\
 &  &  & 5!!x\\
 &  &  &  & \ddots\\
 &  &  &  &  & \left(2N-3\right)!!x\\
  };
\node[xshift=-72pt,yshift=59pt] {$1$};
\node[xshift=-55pt,yshift=59pt] {$2$};
\node[xshift=-31pt,yshift=59pt] {$3$};
\node[xshift=-12pt,yshift=59pt] {$4$};
\node[xshift=12pt,yshift=59pt] {$\cdots$};
\node[xshift=48pt,yshift=59pt] {$N$};
\node[xshift=-95pt,yshift=40pt] {$0$};
\node[xshift=-95pt,yshift=26pt] {$1$};
\node[xshift=-95pt,yshift=11pt] {$2$};
\node[xshift=-95pt,yshift=-4pt] {$3$};
\node[xshift=-95pt,yshift=-18pt] {$\vdots $};
\node[xshift=-102pt,yshift=-38pt] {$N-1$};
\node[xshift=-40pt,yshift=-30pt] {\LARGE$0$};
\end{tikzpicture}
\end{equation*}

From (\ref{eq:11}), we obtain
\begin{align}
 & a_{1}\left(N+1,x\right)\label{eq:16}\\
 & =Na_{0}\left(N,x\right)+xa_{1}\left(N,x\right)\nonumber \\
 & =Na_{0}\left(N,x\right)+x\left(N-1\right)a_{0}\left(N-1,x\right)+x^{2}a_{1}\left(N-1,x\right)\nonumber \\
 & \vdots\nonumber \\
 & =\sum_{i=0}^{N-2}x^{i}\left(N-i\right)a_{0}\left(N-i,x\right)+x^{N-1}a_{1}\left(2,x\right)\nonumber \\
 & =\sum_{i=0}^{N-2}x^{i}\left(N-i\right)a_{0}\left(N-i,x\right)+x^{N-1}x\nonumber \\
 & =\sum_{i=0}^{N-1}x^{i}\left(N-i\right)a_{0}\left(N-i,x\right),\nonumber
\end{align}

\begin{align}
 & a_{2}\left(N+1,x\right)\label{eq:17}\\
 & =\left(N+1\right)a_{1}\left(N,x\right)+xa_{2}\left(N,x\right)\nonumber \\
 & =\left(N+1\right)a_{1}\left(N,x\right)+xNa_{1}\left(N-1,x\right)+x^{2}a_{2}\left(N-1,x\right)\nonumber \\
 & \vdots\nonumber \\
 & =\sum_{i=0}^{N-3}x^{i}\left(N+1-i\right)a_{1}\left(N-i,x\right)+x^{N-2}a_{2}\left(3,x\right)\nonumber \\
 & =\sum_{i=0}^{N-3}x^{i}\left(N+1-i\right)a_{1}\left(N-i,x\right)+3x^{N-2}a_{1}\left(2,x\right)\nonumber \\
 & =\sum_{i=0}^{N-2}x^{i}\left(N+1-i\right)a_{1}\left(N-i,x\right),\nonumber
\end{align}
and
\begin{align}
 & a_{3}\left(N+1,x\right)\label{eq:18}\\
 & =\left(N+2\right)a_{2}\left(N,x\right)+xa_{3}\left(N,x\right)\nonumber \\
 & =\left(N+2\right)a_{2}\left(N,x\right)+x\left(N+1\right)a_{2}\left(N-1,x\right)+x^{2}a_{3}\left(N-1,x\right)\nonumber \\
 & \vdots\nonumber \\
 & =\sum_{i=0}^{N-4}x^{i}\left(N-i+2\right)a_{2}\left(N-i,x\right)+5x^{N-3}a_{2}\left(3,x\right)\nonumber \\
 & =\sum_{i=0}^{N-3}x^{i}\left(N-i+2\right)a_{2}\left(N-i,x\right).\nonumber
\end{align}

Continuing this process, we get
\begin{equation}
a_{j}\left(N+1,x\right)=\sum_{i=0}^{N-j}x^{i}\left(N-i+j-1\right)a_{j-1}\left(N-i,x\right),\label{eq:19}
\end{equation}
where $j=1,2,\dots,N-1$.

Now, we give explicit expressions for $a_{j}\left(N+1,x\right)$ $\left(j=1,2,\dots,N-1\right).$
From (\ref{eq:14}) and (\ref{eq:16}), we can easily derive the following
equation:
\begin{align}
a_{1}\left(N+1,x\right) & =\sum_{i_{1}=0}^{N-1}x^{i_{1}}\left(N-i_{1}\right)a_{0}\left(N-i_{1},x\right)\label{eq:20}\\
 & =x^{N}\sum_{i_{1}=0}^{N-1}\left(N-i_{1}\right).\nonumber
\end{align}

By (\ref{eq:17}), (\ref{eq:18}) and (\ref{eq:19}), we get
\begin{align}
a_{2}\left(N+1,x\right) & =\sum_{i_{2}=0}^{N-2}x^{i_{2}}\left(N-i_{2}+1\right)a_{1}\left(N-i_{2},x\right)\label{eq:21}\\
 & =x^{N-1}\sum_{i_{2}=0}^{N-2}\sum_{i_{1}=0}^{N-2-i_{2}}\left(N-i_{2}+1\right)\left(N-i_{2}-i_{1}-1\right),\nonumber
\end{align}
\begin{align}
a_{3}\left(N+1,x\right) & =\sum_{i_{3}=0}^{N-3}x^{i_{3}}\left(N-i_{3}+2\right)a_{2}\left(N-i_{3},x\right)\label{eq:22}\\
 & =x^{N-2}\sum_{i_{3}=0}^{N-3}\sum_{i_{2}=0}^{N-3-i_{3}}\sum_{i_{1}=0}^{N-3-i_{3}-i_{2}}\left(N-i_{3}+2\right)\left(N-i_{3}-i_{2}\right)\nonumber\\
 &\times \left(N-i_{3}-i_{2}-i_{1}-2\right),\nonumber
\end{align}
and
\begin{align}
a_{4}\left(N+1,x\right) & =\sum_{i_{4}=0}^{N-4}x^{i_{4}}\left(N-i_{4}+3\right)a_{3}\left(N-i_{4},x\right)\label{eq:23}\\
 & =x^{N-3}\sum_{i_{4}=0}^{N-4}\sum_{i_{3}=0}^{N-4-i_{4}}\nonumber\\
 &\times \sum_{i_{2}=0}^{N-4-i_{4}-i_{3}}\sum_{i_{1}=0}^{N-4-i_{4}-i_{3}-i_{2}}\left(N-i_{4}+3\right)\left(N-i_{4}-i_{3}+1\right)\nonumber \\
 & \relphantom =\times\left(N-i_{4}-i_{3}-i_{2}-1\right)\left(N-i_{4}-i_{3}-i_{2}-i_{1}-3\right).\nonumber
\end{align}

Continuing this process, we get
\begin{align}
 & a_{j}\left(N+1,x\right)\label{eq:24}\\
 & =x^{N-j+1}\sum_{i_{j}=0}^{N-j}\sum_{i_{j-1}=0}^{N-j-i_{j}}\cdots\sum_{i_{1}=0}^{N-j-i_{j}-\cdots-i_{2}}\prod_{k=1}^{j}\left(N-i_{j}-\cdots-i_{k}-\left(j-\left(2k-1\right)\right)\right).\nonumber
\end{align}

Therefore, we obtain the following theorem.
\begin{thm}
\label{thm:1} For $N\in\mathbb{N}$, the linear differential equations
\[
F^{\left(N\right)}=\left(\frac{d}{dt}\right)^{N}F\left(t,x\right)=\left(\sum_{i=N}^{2N-1}a_{i-N}\left(N,x\right)\left(1-2t\right)^{-\frac{i}{2}}\right)F
\]
has a solution $F=F\left(t,x\right)=e^{x\left(1-\sqrt{1-2t}\right)},$
where
\begin{align*}
a_{0}\left(N,x\right) & =x^{N},\quad a_{N-1}\left(N,x\right)=\left(2n-3\right)!!x,\\
a_{j}\left(N,x\right) & =x^{N-j}\sum_{i_{j}=0}^{N-j-1}\sum_{i_{j-1}=0}^{N-j-1-i_{j}}\cdots\sum_{i_{1}=0}^{N-j-1-i_{j}-\cdots-i_{2}}\\
 & \relphantom =\times\left(\prod_{k=1}^{j}\left(N-i_{j}-i_{j-1}-\cdots-i_{k}-\left(j-\left(2k-2\right)\right)\right)\right).
\end{align*}

\end{thm}
Recall the the reverse Bessel polynomials $p_{k}\left(x\right)$ are
given by the generating function as
\begin{align}
F & =F\left(t,x\right)=e^{x\left(1-\sqrt{1-2t}\right)}\label{eq:25}\\
 & =\sum_{k=0}^{\infty}p_{k}\left(x\right)\frac{t^{k}}{k!}.\nonumber
\end{align}

Thus, by (\ref{eq:25}), we get
\begin{align}
F^{\left(N\right)} & =\left(\frac{d}{dt}\right)^{N}F\left(t,x\right)\label{eq:26}\\
 & =\sum_{k=N}^{\infty}p_{k}\left(x\right)\left(k\right)_{N}\frac{t^{k-N}}{k!}\nonumber \\
 & =\sum_{k=0}^{\infty}p_{k+N}\left(x\right)\left(k+N\right)_{N}\frac{t^{k}}{\left(k+N\right)!}\nonumber \\
 & =\sum_{k=0}^{\infty}p_{k+N}\left(x\right)\frac{t^{k}}{k!}.\nonumber
\end{align}

On the other hand, by Theorem 1, we get
\begin{align}
F^{\left(N\right)} & =\left(\sum_{i=N}^{2N-1}a_{i-N}\left(N,x\right)\left(1-2t\right)^{-\frac{i}{2}}\right)F\label{eq:27}\\
 & =\sum_{i=N}^{2N-1}a_{i-N}\left(N,x\right)\left(\sum_{l=0}^{\infty}\left(-\frac{i}{2}\right)_{l}\frac{\left(-2t\right)^{l}}{l!}\right)\left(\sum_{m=0}^{\infty}p_{m}\left(x\right)\frac{t^{m}}{m!}\right)\nonumber \\
 & =\sum_{k=0}^{\infty}\left\{ \sum_{i=N}^{2N-1}a_{i-N}\left(N,x\right)\sum_{l=0}^{k}\binom{k}{l}2^{l}\left(\frac{i}{2}+l-1\right)_{l}p_{k-l}\left(x\right)\right\} \frac{t^{k}}{k!}.\nonumber
\end{align}

Therefore, by (\ref{eq:26}) and (\ref{eq:27}), we obtain the following
theorem.
\begin{thm}
\label{thm:2} For $k\in\mathbb{N}\cup\left\{ 0\right\} $, and $N\in\mathbb{N}$,
we have
\[
p_{k+N}\left(x\right)=\sum_{i=N}^{2N-1}a_{i-N}\left(N,x\right)\sum_{l=0}^{k}\binom{k}{l}2^{l}\left(\frac{i}{2}+l-1\right)_{l}p_{k-l}\left(x\right),
\]
where $\left(x\right)_{n}=x\left(x-1\right)\left(x-2\right)\cdots\left(x-n+1\right)$,
$\left(n\ge1\right)$, and $\left(x\right)_{0}=1$. \end{thm}

\bibliographystyle{amsplain}
\nocite{*}

\end{document}